\documentclass[final,3p,times]{elsarticle}
\usepackage{amssymb}
\usepackage{datetime}

\usepackage{listings}
\usepackage{color}

\newcommand{\etal}{\textit{et~al}. }
\newcommand{\etals}{\textit{et~al}.'s}

\definecolor{bg}{rgb}{0.95,0.95,0.95}
\definecolor{mygray}{rgb}{0.5,0.5,0.5}
\definecolor{mygreen}{rgb}{0, 0.6, 0.0}
\definecolor{mymauve}{rgb}{0.58,0,0.82}
\newcounter{bla}

\journal{IJMPC: https://doi.org/10.1142/S0129183121500170}

\begin{document}

\begin{frontmatter}

\title{\emph{AIMpy}: A Python code to solve Schrödinger-like equations with The Asymptotic Iteration Method}

\author{Mesut Karako\c{c}\corref{author}}

\cortext[author] {Corresponding author.\\\textit{E-mail address:} karakoc@akdeniz.edu.tr or mesutkarakoc@gmail.com}
\address{Department of Physics, Faculty of Science, Akdeniz University, TR 07070, Antalya, Turkey}

\begin{abstract}
This paper is dedicated to present an open-source program so-called \emph{AIMpy} built on Python language. \emph{AIMpy} is a solver for Schr\"{o}dinger-like differential equations using Asymptotic Iteration Method (AIM). To confirm the code works seamlessly, it has been shown through the paper with recalculation of some previously studied eigenvalue examples that the code can reproduce their results very well.
\end{abstract}


\end{frontmatter}
%
%
{\bf PROGRAM SUMMARY}
  
\begin{small}
\noindent
{\em Program Title: AIMpy}                                          \\
{\em Licensing provisions: GPLv3 }                                   \\
{\em Programming language: Python 2.x and 3.x}                                   \\
{\em Supplementary material: A Jupyter notebook with examples.}     \\
{\em Nature of problem: Solving Schr\"{o}dinger-like differential equations}\\
{\em Solution method: Asymptotic iteration method}\\
  {\em Source code: https://github.com/mkarakoc/aim}
   \\
\end{small}

\section{Introduction}
\label{sec:intro}

Solving differential equations is of great interest since the invention of calculus by
Newton and Leibniz. After two centuries and fifty years of the invention of calculus, Schr\"{o}dinger introduced his differential equation which now carries his name. It has been applied widely to many physical (quantum) system. One general class of these physical systems are bound systems such as hydrogen atom, molecules, and nucleons in a nuclei, etc. Solutions of such systems mean to obtain their eigenvalues and (or) eigenfunctions.

AIM is one of the various solving methods of the corresponding differential equations for bound systems which is developed by Ciftci, Hall and Saad \cite{Ciftci2003_JPA,Ciftci2005_JPA,Ciftci2005_PRA,Ciftci2005_PLA}. It has many successful application in the literature. Depending on the nature of a differential equation AIM can give analytical or numerical solutions.

Schr\"{o}dinger, Dirac, Klein-Gordon, or Duffin-Kemmer-Petiau (DKP) all these differential wave equations and any others have to take second-order homogeneous linear differential equations of the form 
\begin{equation}
\label{eq:aim}
y'' = \lambda_0(x) y' + s_0(x) y,
\end{equation}
as it is mentioned in Ciftci \emph{et al} \cite{Ciftci2003_JPA} to be solvable by AIM.

In Ref. \cite{Boztosun2006_JMP,Yasuk2006_JMP}, one can find a successful application of AIM to DKP equations which gives analytical solutions for harmonic oscillator and Coulomb potential cases while anharmonic potential case needs a perturbation treatment in the frame work of AIM. Another study which is performed on Schr\"{o}dinger Equation for Makarov Potential \cite{Bayrak2008_IJTP} resulted in analytical solutions. Also there is an approximate analytical solution of Dirac eqaution with Hulthen potential as presented in \cite{Soylu2007_JMP}. There is also a study \cite{Boztosun2007_CPL} which shows AIM can be simplified for a group of differential equations and analytical solutions can be obtained almost immediately.

Although AIM can solve many differential equations analytically in many other cases this might not be possible due to the nature of that particular differential equation, but in these cases, AIM might offer numerical solutions. Cases of numerically obtained eigenvalue can be found in a Schr\"{o}dinger equation study with a Yukawa potential \cite{Karakoc2006_IJMPE}, a DKP equation study which was performed on sextic oscillator potential \cite{Yasuk_2008}, and the references therein. An interesting numerical application of AIM has been performed to Kerr-(A)dS black holes using Kerr-(A)dS angular separation equation \cite{Cho2009_PRD}.

Almost all examples given above (to the best of my knowledge) were investigated using the codes written either in Mathematica or Maple or some similar closed source software. But, \emph{AIMpy} \cite{aimpy} is an open-source \cite{opensource} code based on Python language \cite{Python}, the code is as fast as an unpublished alternative code which has been written in Mathematica by the author of this study.

\emph{AIMpy} \cite{aimpy} code is designed to solve numerical cases and it relies on the improved version of AIM which is developed by Cho \etal \cite{Cho2010_CQG}. The following sections are dedicated to short summary of the description of AIM and improved AIM, and some examples of recalculations using \emph{AIMpy} \cite{aimpy} for the problems chosen from the literature \cite{Karakoc2006_IJMPE,Yasuk_2008,Cho2010_CQG,Bayrak2007_IJQC}.

\section{Model}
\label{sec:model}

\subsection{AIM}
\label{sec:aim}

If one focuses on solving a 1-dimensional Schr\"{o}dinger equation with a smooth changing potential,
this equation can be represented in the following form,
\begin{equation}
\label{eq:1dSch}
\frac{d^2\psi(x)}{dx} + V(x)\psi(x) = E \psi(x)~.
\end{equation}
To be able to perform AIM on Eq. \ref{eq:1dSch} it needs to be transformed into form of Eq. \ref{eq:aim}. This can be achieved by proposing the following relation
\begin{equation}
\label{eq:aimFunc}
\psi(x) = f(x) y(x).
\end{equation}
This relation produces an equation when inserted into Eq. \ref{eq:1dSch} in the form of Eq.~\ref{eq:aim} where
\begin{eqnarray}
\label{eq:l0s0}
\lambda_0(x) &=& - \frac{2 \frac{d}{d x} f{\left (x \right )}}{f{\left (x \right )}} \texttt{, and } \nonumber \\
s_0(x)       &=& E - V{\left (x \right )} - \frac{\frac{d^{2}}{d x^{2}}  f{\left (x \right )}}{f{\left (x \right )}}.
\end{eqnarray}
The process of obtaining $\lambda_0$ and $s_0$ is very similar to  1D~Schr\"{o}dinger case for any differential equation once they are transformed into form of Eq.~\ref{eq:aim}. Proposition of $f(x)$ function is important to obtain the solution, usually one should analyze the asymptotic behaviors of the differential equation: the examples of how $f(x)$ is proposed in can be seen in Refs. \cite{Ciftci2003_JPA,Ciftci2005_JPA,Ciftci2005_PLA,Ciftci2005_PRA,Bayrak2008_IJTP,Karakoc2006_IJMPE,Yasuk_2008,Bayrak2007_IJQC}.

After this point, iterative part of AIM starts due to one has to take derivatives of the Eq.~\ref{eq:aim} and put previously obtained lower degree derivatives in the new higher degree differential equation. This leads to iterative equations to $\lambda_k$ and $s_k$:
\begin{eqnarray} 
\label{eq:lksk}
\lambda _ {k} (x) & = &
\lambda _ {k - 1} ^ {\prime} (x) + s _ {k - 1} (x) + \lambda _ {0} (x) \lambda _ {k - 1} (x) \\ 
s _ {k} (x) & = &
s _ {k - 1} ^ {\prime} (x) + s _ {0} (x) \lambda _ {k - 1} (x) , \quad \textrm{where } k = 1,2,3 , \dots \nonumber
\end{eqnarray}
Asymptotic behavior of the method comes into play when $k$ is large enough and, it gives the following relation,
\begin{equation}
\label{eq:alpha_x}
\frac { s _ { k } (x) } { \lambda _ { k } (x) } = \frac { s _ { k - 1 } (x) } { \lambda _ { k - 1 } (x) } \equiv \alpha (x) .
\end{equation}
A new form of this equation is defined as $\delta_k$,
\begin{equation}
\label{eq:quantization}
\delta _ { k }  \equiv \lambda _ { k }(x)  s _ { k - 1 } ( x ) - \lambda _ { k - 1 } (x) s _ { k } (x),
\end{equation}
and $\delta_k=0$ is called as `quantization condition' \cite{Cho2010_CQG}. Roots of the quantization condition give eigenvalues of any differential equation which obeys AIM conditions. Finally, the following integral which details are given in Ref. \cite{Ciftci2003_JPA} gives the general solution of Eq.~\ref{eq:aim},
\begin{eqnarray}
\label{eq:aim_gen_sol}
y(x) &=& \exp \left( - \int ^ {x} \alpha \left(x^ { \prime } \right) dx^ { \prime } \right) \times \nonumber \\ 
&&\left [ C _ { 2 } + C _ { 1 } \int ^ {x} \exp \left( \int ^ {x^ { \prime }} 
\left[ \lambda _ { 0 } \left(x^ { \prime \prime } \right) + 2 \alpha \left(x^ { \prime \prime } \right) \right] dx^ { \prime \prime } \right) dx^ { \prime } \right].
\end{eqnarray}

\subsection{Improved AIM}
\label{sec:imp_aim}

An apparent `weakness' of AIM is its iterative nature, because one has to take derivatives of $\lambda_{k-1}$ and $s_{k-1}$ to obtain $\lambda_k$ and $s_k$. These derivatives are added to every new iteration and for sufficiently large $k$'s $\lambda_k$ and $s_k$ become huge mathematical terms which might be represented by with very a big order polynomial. This might brings difficulties depending on the case for both numerical and symbolical calculations. The part for the symbolic difficulty would be time and memory consuming, and the part for the numerical difficulty would be precision and accuracy loss.

To overcome this issues Cho \etal \cite{Cho2010_CQG} developed an `improved' version of AIM. According to the improved version, $\lambda_k$ and $s_k$ can be written as Taylor series,
\begin{eqnarray} 
\lambda _ {n} ( x ) &= & \sum _ { i = 0 } ^ { \infty } c _ {n} ^ { i } ( x - x_0 ) ^ { i }, \\ 
s _ {n} ( x ) &= & \sum _ { i = 0 } ^ { \infty } d _ {n} ^ { i } ( x - x_0 ) ^ { i }, 
\end{eqnarray}
where $c _ {n} ^ { i }$ and $d _ {n} ^ { i }$ are Taylor coefficients. Using these series in Eqs.~\ref{eq:lksk} gives following recursion relations for the coefficients:
\begin{eqnarray}
c _ {n} ^ { i } &=& ( i + 1 ) c _ {n-1} ^ { i + 1 } + d _ {n-1} ^ { i } + \sum _ { j = 0 } ^ { i } c _ { 0 } ^ { j } c _ {n-1} ^ { i - j },\\
d _ {n} ^ { i } &=& ( i + 1 ) d _ {n-1} ^ { i + 1 } + \sum _ { j = 0 } ^ { i } d _ { 0 } ^ { j } c _ {n-1} ^ { i - j }.
\end{eqnarray}
The quantization condition in Eq. \ref{eq:quantization} can be rewritten with these new recursion relations in the following form,
\begin{equation}
\label{eq:newquantization}
d_ {n} ^ { 0 } c_{n-1} ^ { 0 } - d _ {n-1} ^ { 0 } c _ {n} ^ { 0 } = 0
\end{equation}
and there is no need take derivatives when using the new condition other than the first iteration. Therefore, only $\lambda_0$ and $s_0$ derivatives are taken to obtain $c _ {0} ^ { i }$ and $d _ {0} ^ { i }$ coefficients which are necessary to start the new recursion relation calculations.

\section{The Code}
\label{sec:thecode}

\emph{AIMpy} \cite{aimpy} is designed in an open-source environment (Python language), and also I believe open-source \cite{opensource} philosophy is one of the driving force to the progress of science. Therefore, \emph{AIMpy} \cite{aimpy} is an open-source code with GPL \cite{GPL} license which allows anyone can modify and develop on top of it.

It is designed with two parts: the one including all necessary libraries and Python functions and, the other one includes only AIM calculations. The first part is a Python program called \emph{`asymptotic.py'} and the other one a Jupyter notebook \cite{Jupyter} which one can give any name. However, it is preferred as \emph{`AIMpy\_[name of the problem].ipynb'} for the examples of this paper. The code, the examples and,  installation and usage manuals can be found at https://github.com/mkarakoc/aim. 

\emph{`asymptotic.py'}  has Python libraries which are named IPython \cite{IPython}, \emph{SymPy} (sympy) \cite{SymPy}, \emph{symengine} \cite{SymEngine} and Python-FLINT (flint) \cite{PythonFlint}.
IPython is necessary to show mathematical terms as user-friendly. For example, a parameter named as `beta' will be shown in the outputs as $\beta$. 
\emph{SymPy} and \emph{symengine} are both used to make symbolic calculations, mainly for derivatives of $\lambda_0$ and $s_0$. Both \emph{SymPy} and \emph{symengine} can calculate arbitrary precision numbers and that makes them right tools for AIM since the arbitrary precision (or high precision) calculation is one of the key-point of AIM. \emph{SymPy} and \emph{symengine} both are capable of taking derivatives but \emph{symengine} is significantly faster than \emph{SymPy}. Therefore, \emph{symengine} is preferred for derivatives while \emph{SymPy} is mostly used to convert symbolic variables into \LaTeX~to show outputs as user-friendly using IPython tools. 

Once, $c _ {k} ^ { i }$ and $d _ {k} ^ { i }$ coefficients are obtained using these two libraries then the quantization condition in Eq. \ref{eq:newquantization} is obtained in an polynomial form depending on $x$ (equation variable) and $E$ (eigenvalue variable) parameters. The roots of this polynomial are calculated using \emph{Python-FLINT} or shortly flint. As it is mentioned in the reference page \cite{PythonFlint}: \emph{Python-FLINT} is an \emph{`Python extension module wrapping FLINT (Fast Library for Number Theory) and Arb (arbitrary-precision ball arithmetic).'} The roots can have arbitrary precision since both \emph{FLINT} \cite{flint} and \emph{Arb} \cite{arb} are arbitrary precision libraries. Therefore, it can be claimed that the eigenvalues that would be calculated with \emph{AIMpy} \cite{aimpy} will be very reliable. 

Second part of the code is in a Jupyter notebook \cite{Jupyter} which actually is an environment to create an user-friendly interface and to just take care of the physics problem interested in rather than the code itself.

\section{Examples}
\label{sec:examples}

The problems solvable by AIM can be classified as analytically and numerically solvable ones, but this version of \emph{AIMpy} does not have analytical solver property, yet. Therefore, the below examples are selected from the literature \cite{Karakoc2006_IJMPE,Bayrak2007_IJQC,Yasuk_2008,Cho2010_CQG} to present \emph{AIMpy} for numerical applications of AIM. The full details of the examples are not presented in this study. One shall see the details of the equations, the potentials and the definitions of their parameters in the reference papers since these examples are focused on to confirm \emph{AIMpy} reproduce results of earlier studies very well.

The first example is given in a more detailed fashion to show how the \emph{AIMpy} code would look like and how it is used to solve the eigenvalue problems. The other examples are kept more abstract since they are just proof of the success of \emph{AIMpy}.

\subsection{The Yukawa Potential}
\label{sec:yukawa}

\begin{equation}
\label{eq:Yukawa}
V(r)=-\frac{A}{r} \exp (-\alpha r)
\end{equation}

The first example is the solution of the radial Schr\"{o}dinger equation (Eq.~\ref{eq:RadSch}) for a Yukawa type potential (Eq.~\ref{eq:Yukawa}) which was studied in Ref. \cite{Karakoc2006_IJMPE} using standard AIM through the unpublished Mathematica code written by the author. The eigenvalues solution of the system is given below such as seen in a Jupyter notebook. One can get an idea of how the other examples also have been solved with \emph{AIMpy} \cite{aimpy}.

\begin{equation}
\label{eq:RadSch}
\frac{d^{2} R_{n}(r)}{d r^{2}}+\frac{2 m}{\hbar^{2}}\left(E_{nL}-\frac{L(L+1) \hbar^{2}}{2 m r^{2}} - V(r)\right) R_{n}(r)=0
\end{equation}
\begin{lstlisting}
# Python program to use AIM tools
from asymptotic import *

# symengine (symbolic) variables 
# for lambda_0 and s_0 
En,m,hbar,L,r,r0 = se.symbols("En,m,hbar,L,r,r0")
beta,alpha,A = se.symbols("beta,alpha,A")
\end{lstlisting}

In the first step, one should import \emph{``asymptotic.py"}, then define the necessary paramaters as symbolic variables of \emph{symengine}. Next step is the writing of $\lambda_0$ and $s_0$ (Eq. \ref{eq:yukawa_l0s0}) with these variables as seen below. It should be noted here that $\lambda_0$ and $s_0$ must be obtained by user as it has been explained in the section \ref{sec:aim}.
\begin{eqnarray}
\label{eq:yukawa_l0s0}
\lambda_0 &=& 2 \beta - \frac{2}{r} \nonumber \\
s_0 &=&  - \frac{2m }{\hbar^{2}}\left(E_{nL}  -\frac{A e^{- \alpha r}}{r}\right) + \frac{L(L+1)}{r^{2}}  + \frac{2 \beta}{r} - \beta^{2} 
\end{eqnarray}
\begin{lstlisting}
# 1st step:
# lambda_0 and s_0
l0 = 2*beta - 2/r
s0 = -2*m/hbar**2*(En-A*se.exp(-alpha*r)/r) + L*(L+1)/r**2 + (2*beta)/r - beta**2 
\end{lstlisting}

The user should be aware of all numbers are in l0~$\equiv\lambda_0$ and s0~$ \equiv s_0$ are in infinite precision until the iterative calculation process starts. It should be avoided to put numerical values of the symbolic parameters directly into the l0 and s0. Therefore it would be possible to use symbolic form of l0 and s0 later on if it is needed to make  calculations for different numerical values of the parameters. 
\begin{lstlisting}
# 2nd step:
## Case: A = 4, L=0
# values of variables
nA = o*4
nL = o* 0
nalpha = o* 2/10 
nhbar, nm = o* 1, o*1/2

nbeta = o* 3
nr0 = o* 1/nbeta

# parameters of lambda_0 (pl0) and s_0 (ps0)
pl0 = {beta: nbeta}
ps0 = {beta: nbeta, alpha: nalpha, A:nA, m: nm, L: nL, hbar: nhbar, r0: nr0}
\end{lstlisting}

Third step is the preparation of the numerical values of the parameters.
It is easy to understand that the variable names that start with \emph{``n"} are for the numerical values. These values are connected with the symbolic variables through pl0 and ps0 which are Python dictionaries. The benefit of this usage l0 and s0 still have been kept in their symbolic forms. As it is mentioned earlier all numerical values will be kept in infinite precision to assure this an object named as \emph{``o"} is created. The value of is {\bf 1} with infinite precision and it has an \emph{symengine} integer type.  For example, \emph{o * 1/2} will provide that \emph{1/2} will always be stay as it is till the precision of it to changed an finite number. 

The existence of parameters $A, L, \alpha, \hbar, m$ are obvious and their numerical values are taken from Ref. \cite{Karakoc2006_IJMPE}. The $\beta$ parameter comes from the asymptotic form of the proposed wave-function (Ref. \cite{Karakoc2006_IJMPE}) in Eq. \ref{eq:yukawa_wave}. Its value is just an arbitrary number and usually chosen to make iterative process to converge as early as possible to eigenvalues. $r_0 = 1/\beta$ is chosen as the maximum of the asymptotic part of the wave-function as it is suggested in Ref. \cite{Karakoc2006_IJMPE}. 
\begin{equation}
\label{eq:yukawa_wave}
R_{n}(r)=r \exp (-\beta r) f(r)
\end{equation}

\begin{lstlisting}
# 3rd step:
# pass lambda_0, s_0 and variable values to aim class
yukawa_A4L0 = aim(l0, s0, pl0, ps0)
yukawa_A4L0.display_parameters()
yukawa_A4L0.display_l0s0(0)
yukawa_A4L0.parameters(En, r, nr0, nmax=201, nstep=10, dprec=500, tol=1e-101)
\end{lstlisting}

The fourth step is to create a \emph{Python object} with using a \emph{Python class} named as \emph{aim}. This class takes four inputs l0, s0 are symbolic representations of $\lambda_0$ and $s_0$ in Python language and, pl0, ps0 are Python dictionaries they contain numerical values for l0 and s0, respectively. The numerical values passed to l0 and s0 inside of the \emph{aim} class. In this step, the Python object created by \emph{aim} is assigned to \emph{``yukawa\_A4L0"} variable. These name can be any arbitrary name which obeys Python naming rules. But in this example, \emph{``yukawa"} stands for the potential name and, the rest of name gives an idea of which case is user studying on. \emph{``yukawa\_A4L0.display\_parameters()"} and \emph{``yukawa\_A4L0.display\_l0s0(0)"} lines are not necessary for the solution, they are to show the numerical values of the parameters and $\lambda_0$ and $s_0$ in much more a human readable representation. The last line contains the name of the important variables. The defitions of the variables are given in Table \ref{tab:var1}.

\begin{table}[!hbtp]
	\begin{center}
		\caption{Names of the variables can be any Pythonic variable name, but it can also be chosen similar to the variables of a differential equations in interest.}
		\label{tab:var1}
		\begin{tabular}{|l|l|}
			\hline 
			Variable & Definitions\\ 
			Name&  \\ 
			\hline 
			En	& eigenvalue of the differential equation\\ 
			r	& the differential equation variable which is used for derivatives\\ 
			nr0	& a particular value for r, usually minimum of the potential or \\
			& the maximum of the asymptotic part of the wave-function \\ 
			nmax& the maximum iteration number\\
			nstep& skip every iteration with the value given to nstep\\
			dprec& the finite precision of all numerical values\\
			tol & tolerance for convergence to eigenvalues \\
			\hline 
		\end{tabular} 
	\end{center}
\end{table}
\begin{lstlisting}
# 4th step:
# create coefficients for improved AIM
yukawa_A4L0.c0()
yukawa_A4L0.d0()
yukawa_A4L0.cndn()
\end{lstlisting}

The fifth step is the calculations of Taylor series coefficients according to improved version of AIM (Ref. \cite{Cho2010_CQG}). $c _ {0} ^ { i }$ and $d _ {0} ^ { i }$ and are calculated with the lines \emph{``yukawa\_A4L0.c0()"} and \emph{``yukawa\_A4L0.d0()"}, respectively. $c _ {n} ^ { i }$ and $d _ {n} ^ { i }$ (where $n = 1, 2, 3, \ldots$) are calculated with the last line (\emph{``yukawa\_A4L0.cndn()"}).

The last step is the following code line. As it is mentioned earlier, \emph{``Arb"} library is used through \emph{``Python-FLINT"} to obtain the eigenvalues. Therefore, the method to obtain the eigenvalues named as \emph{``get\_arb\_root"}. In this example, eigenvalues are real negative numbers where their fraction part has 20 digits. There might be other roots like real positive or complex, but they are filtered with \emph{``showRoots='-r'"}. The size of the fraction part is defined with \emph{``printFormat"}, and it is obvious \emph{.20f} in the example stands for 20 digits. 

\begin{lstlisting}
# 5th step:
yukawa_A4L0.get_arb_roots(showRoots='-r', printFormat="{:25.20f}")
\end{lstlisting}

Finally, the results of the calculation are given below. It can be seen from the outputs that the fast and successful convergence of the iterations. In the earlier study that I have contributed \cite{Karakoc2006_IJMPE} lesser fractional digits are presented which it was enough for the goal of that study (\emph{see Table \ref{tab:Yukawa}}), but it is aimed to show the how successful AIMpy \cite{aimpy} in the present one.

\lstset{
	basicstyle=\footnotesize, 
	backgroundcolor=\color{bg},
	commentstyle=\color{mygreen},
	stringstyle=\color{mymauve},
	keywordstyle=\color{blue},
	breaklines=false,
	language=Python,
	frame=single
}
\begin{lstlisting}
# 6th (last) step:
iteration          E00                        E10                         E20
   001   -2.22608382037941285277                                                      
   021   -3.25646424490443058235   -0.39503930505322610330                            
   041   -3.25646424490722525404   -0.39942480823044205382   -0.00752426627166027709
   061   -3.25646424490722525404   -0.39942617037004934161   -0.02379348447266779420
   081   -3.25646424490722525404   -0.39942617065529516410   -0.02544763075010264581
   101   -3.25646424490722525404   -0.39942617065535110306   -0.02563879250512245213
   121   -3.25646424490722525404   -0.39942617065535111388   -0.02566135365381434051
   141   -3.25646424490722525404   -0.39942617065535111388   -0.02566401769819825872
   161   -3.25646424490722525404   -0.39942617065535111388   -0.02566433195250032477
   181   -3.25646424490722525404   -0.39942617065535111388   -0.02566436900553236619
   201   -3.25646424490722525404   -0.39942617065535111388   -0.02566437337375074847
CPU times: user 21.1 s, sys: 24.8 ms, total: 21.2 s, Wall time: 21.2 s
\end{lstlisting}

\begin{table}[!htbp]
	\caption{The energy eigenvalues ($E_{nL}$) of the Yukawa potential (Eq. \ref{eq:Yukawa}) where $\hbar=2m=1$,$\alpha=0.2 fm^{-1}$, $\beta=3$ and  $n=0$ for all $L$ values. The results of the Refs. \cite{Gonul2006_PS,Chakrabarti2001_PLA,Chakrabarti2001_PLA} are just given to show the results of the other methods. The result of Ref. \cite{Karakoc2006_IJMPE} is the one that one should really compare of the present study.}
	\label{tab:Yukawa}
	\begin{center}
		\begin{tabular}{rrrrrrr}
			\hline
			$A$ &$L$& Present study ($E_{nL}$) & AIM\cite{Karakoc2006_IJMPE} & SUSY \cite{Gonul2006_PS}& Numerical\cite{Chakrabarti2001_PLA} & Analytical\cite{Chakrabarti2001_PLA}\\ 
			\hline \vspace{4pt}
			4 & 0 &   -3.25646424490722525404 &  -3.256464 &   -3.2563 &   -3.2565 &   -3.2199 \\
			8 & 0 &   14.45812571278417740340 & -14.458126 &  -14.4581 &  -14.4571 &  -14.4199 \\ \vspace{4pt}
			& 1 &   -2.58369238520910751079 &  -2.583692 &   -2.5830 &   -2.5836 &   -2.4332 \\
			16 & 0 &  -60.85903282302551371170 & -60.859033 &  -60.8590 &  -60.8590 &  -60.8193 \\ \vspace{4pt}
			& 1 &  -12.99103533706039481539 & -12.991035 &  -12.9908 &  -12.9910 &  -12.8375 \\
			24 & 0 & -139.25934814287272696744 &-139.259348 & -139.2590 & -139.2594 & -139.2201 \\
			& 1 &  -31.39381360113816006395 & -31.393814 &  -31.3937 &  -31.3938 &  -31.2385 \\
			& 2 &  -11.59594912057730331414 & -11.595949 &  -11.5951 &  -11.5959 &  -11.2456 \\
			\hline
		\end{tabular}
	\end{center}
\end{table}

\subsection{The Exponential Cosine Screened Coulomb (ECSC) Potential}
\label{sec:exppot}

\begin{equation}
\label{eq:ecsc}
V(r)=-\frac{A}{r} e^{-\delta r} \cos (\delta r)
\end{equation}

The ECSC potential in Eq. \ref{eq:ecsc} for the radial Schr\"{o}dinger equation (Eq.~\ref{eq:RadSch}) has been analyzed in Refs. \cite{Bayrak2007_IJQC,Ikhdair1993_ZPD,Meyer1985_JPA,Ikhdair2007_JMC} with different solution methods. Ref. \cite{Bayrak2007_IJQC} is followed to redo the calculations with AIMpy \cite{aimpy}. Therefore, one can find in Bayrak \etal \cite{Bayrak2007_IJQC} how $\lambda_{ 0 }$, $s_0$ (Eqs. \ref{eq:ecsc_l0s0}): 
\begin{eqnarray}
\label{eq:ecsc_l0s0}
\lambda_0 &=& 2 \beta - \frac{2 \left(L + 1\right)}{r},\\ 
s_0 &=& - \frac{2m }{\hbar^{2}}E_{nL} - \beta^{2} +  \frac{2 \beta (L+1)}{r} -\frac{A_{1}}{r}+ A_{2} - A_{3} r^{2} + A_{4} r^{3} - A_{5} r^{4} + A_{6} r^{6},  \nonumber
\end{eqnarray}
where,
\begin{equation}
A_{1}=\frac{2 m}{\hbar^{2}} A, \quad 
A_{2}=A_{1} \delta, \quad 
A_{3}=\frac{A_{1}}{3}\delta^{3}, \quad 
A_{4}=\frac{A_{1}}{6}\delta^{4}, \quad 
A_{5}=\frac{A_{1}}{30}\delta^{5}, \quad 
A_{6}=\frac{A_{1}}{630}\delta^{7}
\end{equation}
and, the terms in them are obtained. 

The results of present study using AIMpy \cite{aimpy} are presented in Table \ref{tab:ExpCos}. The results of the reference paper Ref. \cite{Bayrak2007_IJQC} are not presented here since the present results are almost exactly the same with them other than one or two last digits of the fractional parts of the eigenvalues.

\begin{table}[!htbp]
	\caption{The energy eigenvalues ($E_{nL}$) of the ECSC potential (Eq. \ref{eq:ecsc}) where $A, \hbar=m=1$, $\beta=6/10$. Compare the results with Ref. \cite{Bayrak2007_IJQC}.}
	\label{tab:ExpCos}
	\begin{center}
		\begin{tabular}{cccccc}
			\hline
			$\delta$ & $E_{nL}$ & $E_{1L}$ & $E_{2L}$ & $E_{3L}$ & $E_{4L}$ \\
			\hline
			0.01 &   $s$   & -0.49000099 & -0.11501346 & -0.04561908 & -0.02143746 \\
			&   $p$   &             & -0.11500966 & -0.04561104 & -0.02142437 \\
			&   $d$   &             &             & -0.04559484 & -0.02139798 \\ \vspace{4pt}
			&   $f$   &             &             &             & -0.02135784 \\
			0.02 &   $s$   & -0.48000780 & -0.10510359 & -0.03602510 & -0.01257152 \\
			&   $p$   &             & -0.10507464 & -0.03596760 & -0.01248554 \\
			&   $d$   &             &             & -0.03585066 & -0.01231013 \\ \vspace{4pt}
			&   $f$   &             &             &             & -0.01203814 \\
			0.06 &   $s$   & -0.44020051 & -0.06742086 & -0.00534922 &             \\
			&   $p$   &             & -0.06677740 & -0.00438279 &             \\ \vspace{4pt}
			&   $d$   &             &             & -0.00226187 &             \\
			0.10 &   $s$   & -0.40088477 & -0.03491583 &             &             \\
			&   $p$   &             & -0.03245501 &             &             \\
			\hline
		\end{tabular}
	\end{center}
\end{table}

\subsection{The Sextic Oscillator Potential}
\label{sec:sectpot}

The following homogeneous second-order differential equation for the DKP anharmonic (sextic) oscillator:
\begin{equation}
\left[\frac{\mathrm{d}^{2}}{\mathrm{d} r^{2}}-\frac{J(J+1)}{r^{2}}+\frac{3 m w}{\hbar}-\left(\frac{m^{2} w^{2}}{\hbar^{2}}+\frac{5 q}{\hbar}\right) r^{2}
+\frac{2 m w q}{\hbar^{2}} r^{4}-\frac{q^{2}}{\hbar^{2}} r^{6} \right] F(r)=\frac{1}{\hbar^{2} c^{2}}\left(m^{2} c^{4}-E^{2}\right) F(r),
\end{equation}
can be obtained when the process in Ref. \cite{Yasuk_2008} is followed. This is a second order differential equation which it has a very similar mathematical form to the radial Schr\"{o}dinger equation (Eq.~\ref{eq:RadSch}). Therefore, after this point it is easy to convert this equation following Ref. \cite{Yasuk_2008} and by proposing an asymptotic form for $F(r)$ function to obtain a solvable equation by AIM through AIMpy \cite{aimpy}. Through this way one can obtain following $\lambda_0$ and $s_0$ forms:
\begin{eqnarray}
\lambda_0 &=& 2 \beta_{1} + 4 \beta_{2} r - \frac{2}{r},\\
s_{0}&=&-\beta_{1}^{2}+6 \beta_{2}+\frac{2 \beta 1}{r}-4 \beta_{1} \beta_{2} r-4 \beta_{1}^{2} r^{2}-\frac{1}{\hbar^{2} c^{2}}\left(E_{nJ}^{2}+A_{0}-\frac{A_{1}}{r^{2}}-A_{2} r^{2}+A_{3} r^{4}-A_{4} r^{6}\right), \nonumber
\end{eqnarray}
where
\begin{equation}
A_{0}=m c^{2}\left(3 \hbar w-m c^{2}\right), \quad
A_{1}=\hbar^{2} c^{2} J(J+1), \quad
A_{2}=c^{2}\left(m^{2} w^{2}+5, \hbar q\right), \quad
A_{3}=2 m c^{2} q w, \quad
A_{4}=q^{2} c^{2}.
\end{equation}

The eigenvalues obtained in the present study are presented in Table \ref{tab:Sectic} which are exactly the same with the ones in the reference \cite{Yasuk_2008} paper.

\begin{table}[!htbp]
	\caption{Ground and excited state energies of the sextic oscillator where $\hbar = c = m = 1$ and $q = w = 0.1$. Compare the results with Ref. \cite{Yasuk_2008}.}
	\label{tab:Sectic}
	\begin{center}
		\begin{tabular}{ccccc}
			\hline
			$E_{nJ}$ & $n=0$ & $n=1$ & $n=2$ & $n=3$ \\
			\hline
			$s$ & 1.72356712431419 & 2.59853197666838 & 3.39500190327295 & 4.14085218155052 \\
			$p$ & 2.15692634351980 & 2.98966248683064 & 3.76218401610683 & 4.49009342299632 \\
			$d$ & 2.54665944322791 & 3.35618211593500 & 4.11089275949862 & 4.82457020376795 \\
			$f$ & 2.91026031614511 & 3.70392637878655 & 4.44477489149058 & 5.14682205155551 \\
			\hline
		\end{tabular}
	\end{center}
\end{table}

\subsection{Black Hole Quasinormal Modes}
\label{sec:blackhole}

In this example, the results of the reference paper Ref. \cite{Cho2010_CQG} for black hole quasinormal modes (QNMs) has been reproduced. Improved AIM is one of the important parts of the AIMpy \cite{aimpy} is also presented by Cho \etal \cite{Cho2010_CQG} in the same paper. Cho \etal \cite{Cho2010_CQG} focuses on a field equation of the form
\begin{equation}
\frac{\mathrm{d}^{2} \psi(x)}{\mathrm{d} x^{2}}+\left[\omega^{2}-V(x)\right] \psi(x)=0,
\end{equation}
where V(x) is a master potential of the form \cite{Cho2010_CQG,Berti2009_CQG}
\begin{equation} 
V(r)=f(r)\left[\frac{\ell(\ell+1)}{r^{2}}+\left(1-s^{2}\right)\left(\frac{2 M}{r^{3}}-\frac{\left(4-s^{2}\right) \Lambda}{6}\right)\right],
\end{equation}
and $\mathrm{d} x=\mathrm{d} r / f(r)$ where
\begin{equation} 
f(r)=1-\frac{2 M}{r}-\frac{\Lambda}{3} r^{2}.
\end{equation}

$\lambda_0(\xi)$ and $s_0(\xi)$ can be obtained with the $\xi=1 / r$ variable change \cite{Cho2010_CQG,Moss2002_CQG} and following the process in Cho~\etals~\cite{Cho2010_CQG} study as

\begin{eqnarray}
\lambda_0(\xi) &=& -{1\over p} \left[ p' - {2 i\omega_{n\ell} \over \kappa_1(\xi - \xi_1)}  -2 i \omega_{n\ell} \right], \\
s_0(\xi) &=& {1\over p} \left[ \ell(\ell +1) + 
(1-s^2) \left( 2M\xi - (4-s^2){\Lambda \over 6\xi^2} \right) + 
{i\omega_{n\ell} \over \kappa_1(\xi - \xi_1)^2} \left( {i\omega_{n\ell} \over \kappa_1} +1 \right) + (p' - 2i\omega_{n\ell}) {i\omega_{n\ell} \over \kappa_1(\xi - \xi_1)}  \right].  \nonumber
\end{eqnarray}

The QNMs ($\omega_{n\ell}$) presented in Table \ref{tab:Blackhole} are belongs to present study. Because the results are almost the same with the reference paper Ref. \cite{Cho2010_CQG} as it happened with the previous examples above.

\begin{table}[!htbp]
	\caption{All QNMs ($\omega_{n\ell}$) with $s=2$ in this table obtained for $\ell=2, 3$ and various $\Lambda$ cosmological constants. Compare the results with Ref. \cite{Cho2010_CQG}.}
	\label{tab:Blackhole}
	\begin{center}
		\begin{tabular}{ccccc}
			\hline
			$\ell=2$ & $\Lambda$ & n=1 & n=2 & n=3 \\
			\hline
			& 0.00 & 0.3736717 - 0.0889623i & 0.3467110 - 0.2739149i & 0.3010535 - 0.4782770i \\
			& 0.02 & 0.3383914 - 0.0817564i & 0.3187587 - 0.2491966i & 0.2827322 - 0.4294841i \\
			& 0.04 & 0.2988947 - 0.0732967i & 0.2858409 - 0.2217241i & 0.2599919 - 0.3770922i \\
			& 0.06 & 0.2532892 - 0.0630425i & 0.2457420 - 0.1897910i & 0.2300764 - 0.3191573i \\
			& 0.08 & 0.1974823 - 0.0498773i & 0.1941148 - 0.1497866i & 0.1871198 - 0.2502570i \\
			& 0.09 & 0.1626104 - 0.0413665i & 0.1607886 - 0.1241522i & 0.1570423 - 0.2071172i \\
			& 0.10 & 0.1179164 - 0.0302105i & 0.1172432 - 0.0906409i & 0.1158764 - 0.1511018i \\ \vspace{4pt}
			& 0.11 & 0.0372699 - 0.0096157i & 0.0372493 - 0.0288470i & 0.0372081 - 0.0480784i \\
			$\ell=3$ & $\Lambda$ & n=1 & n=2 & n=3 \\
			\hline
			& 0.00 & 0.5994433 - 0.0927030i & 0.5826438 - 0.2812981i & 0.5516849 - 0.4790928i \\
			& 0.02 & 0.5431149 - 0.0844957i & 0.5307443 - 0.2553631i & 0.5070153 - 0.4320588i \\
			& 0.04 & 0.4800575 - 0.0751464i & 0.4716583 - 0.2263948i & 0.4550106 - 0.3807731i \\
			& 0.06 & 0.4071752 - 0.0641396i & 0.4021706 - 0.1928074i & 0.3920528 - 0.3227693i \\
			& 0.08 & 0.3178048 - 0.0503821i & 0.3154946 - 0.1512490i & 0.3108033 - 0.2524505i \\
			& 0.09 & 0.2618425 - 0.0416439i & 0.2605716 - 0.1249688i & 0.2579976 - 0.2084119i \\
			& 0.10 & 0.1899943 - 0.0303145i & 0.1895170 - 0.0909507i & 0.1885554 - 0.1516089i \\
			& 0.11 & 0.0600915 - 0.0096189i & 0.0600766 - 0.0288567i & 0.0600469 - 0.0480945i \\
			\hline
		\end{tabular}
	\end{center}
\end{table}

\section{Conclusions and outlook}
\label{sec:conclusions}

I have presented \emph{AIMpy} code which is an eigenvalue solver for Schr\"{o}dinger-like differential equations using AIM. \emph{AIMpy} is an open-source \cite{opensource} code using the power of fast and high precision open-source symbolic and numeric calculation libraries. The example cases given in the paper proves the reliability of the code. However, some more features can be added to the code in later versions: like eigenfunction calculations, analytical solver and, direct input of any solvable differential equations by AIM to AIMpy to obtain $\lambda_0$ and $s_0$, \emph{etc.} As the last word, I would like to encourage community contributions or writing their versions of the code through forking the repository on GitHub \cite{aimpy}.

\bibliographystyle{plain}       
\bibliography{AIMpy_mkarakoc}

\end{document}